\documentclass[12pt]
{amsart}

\newtheorem{Theorem}{Theorem}[section]
\newtheorem{Proposition}[Theorem]{Proposition}
\newtheorem{Lemma}[Theorem]{Lemma}

\theoremstyle{remark}
\newtheorem{Remark}[Theorem]{Remark}

\numberwithin{equation}{section}

\allowdisplaybreaks

\begin{document}

\title[A multidimensional $_8\psi_8$ summation]
{A multidimensional generalization of Shukla's
$\hbox{${}_{\boldsymbol 8}\boldsymbol\psi_{\boldsymbol 8}$}$ summation}
\author{Michael Schlosser}
\address{Department of Mathematics, The Ohio State University,
231 West 18th Avenue, Columbus, Ohio 43210, USA}
\email{mschloss@math.ohio-state.edu}
\urladdr{http://www.math.ohio-state.edu/\textasciitilde mschloss}
\subjclass
{Primary 33D15; Secondary 33D67.}
\keywords{basic hypergeometric series, $q$-series, bilateral series,
$_8\psi_8$ summation, $_6\phi_5$ summation, $A_r$ series, $U\!(r)$ series.}
\date{March 5, 2001}

\dedicatory{\rm\small Department of Mathematics, The Ohio State University,\\
231 West 18th Avenue, Columbus, Ohio 43210, USA\\
E-mail: \tt mschloss@math.ohio-state.edu\\
URL: \tt http://www.math.ohio-state.edu/\textasciitilde mschloss}

\begin{abstract}
We give an $r$-dimensional generalization of H.~S.~Shukla's
very-well-poised $_8\psi_8$ summation formula.
We work in the setting of multiple basic hypergeometric series
very-well-poised over the root system $A_{r-1}$, or equivalently, the unitary
group $U\!(r)$. Our proof, which is already new in the
one-dimensional case, utilizes an $A_{r-1}$ nonterminating
very-well-poised $_6\phi_5$ summation by S.~C.~Milne,
a partial fraction decomposition, and analytic continuation.
\end{abstract}

\maketitle

\section{Introduction}

In this article, we provide a multidimensional generalization of
H.~S.~Shukla's \cite[Eq.~(4.1)]{shukla} very-well-poised
$_8\psi_8$ summation.
In the classical (one-dimensional) case, it reads as follows
(cf.\ H.~Exton~\cite[Eq.~(3.8.1.2) or (A.28)]{exton}):
\begin{multline}\label{8psi8gl}
{}_8\psi_8\!\left[\begin{matrix}q\sqrt{a},-q\sqrt{a},b,c,aq^2/c,e,f,g\\
\sqrt{a},-\sqrt{a},aq/b,aq/c,c/q,aq/e,aq/f,aq/g\end{matrix}\,
;q,\frac{a^2}{befg}\right]\\
=\left(1-\frac{(1-be/a)(1-bf/a)(1-bg/a)}{(1-bq/c)(1-bc/aq)(1-befg/a^2)}\right)
\frac{(1-c/bq)(1-bc/aq)}{(1-c/aq)(1-c/q)}\\\times
\frac{(q,aq,q/a,aq/be,aq/bf,aq/bg,aq/ef,aq/eg,aq/fg;q)_{\infty}}
{(aq/b,aq/e,aq/f,aq/g,q/b,q/e,q/f,q/g,a^2q/befg;q)_{\infty}},
\end{multline}
provided $|q|<1$ and $|a^2/befg|<1$.

Here, and throughout the article, we are using the following notations.
For a complex number $q$ with $|q|<1$, 
the {\em $q$-shifted factorial} is defined by
\begin{equation*}
(a;q)_{\infty}:=\prod_{j=0}^{\infty}(1-aq^j),
\end{equation*}
and
\begin{equation}\label{qpochdef}
(a;q)_k:=\frac{(a;q)_{\infty}}{(aq^k;q)_{\infty}},\qquad
\text{where $k$ is an integer}.
\end{equation}
Further, for brevity, we are employing the notation
\begin{equation*}
(a_1,\ldots,a_m;q)_k\equiv (a_1;q)_k\dots(a_m;q)_k,
\end{equation*}
where $k$ is an integer or infinity.
Finally, we use
\begin{equation}\label{defhyp}
{}_t\phi_{t-1}\!\left[\begin{matrix}a_1,a_2,\dots,a_t\\
b_1,b_2,\dots,b_{t-1}\end{matrix}\,;q,z\right]:=
\sum _{k=0} ^{\infty}\frac {(a_1,a_2,\dots,a_t;q)_k}
{(q,b_1,\dots,b_{t-1};q)_k}z^k,
\end{equation}
and
\begin{equation}\label{defhypb}
{}_t\psi_t\!\left[\begin{matrix}a_1,a_2,\dots,a_t\\
b_1,b_2,\dots,b_t\end{matrix}\,;q,z\right]:=
\sum _{k=-\infty} ^{\infty}\frac {(a_1,a_2,\dots,a_t;q)_k}
{(b_1,b_2,\dots,b_t;q)_k}z^k,
\end{equation}
to denote the {\em basic hypergeometric $_t\phi_{t-1}$ series}, and
{\em bilateral basic hypergeometric $_t\psi_t$ series}, respectively.
For a survey of classical results in the theory of basic hypergeometric series,
see G.~Gasper and M.~Rahman~\cite{grhyp}.
For applications of basic hypergeometric series to various areas,
including number theory, combinatorics, and physics, see
G.~E.~Andrews~\cite{andappl,qandrews}.

The $_8\psi_8$ summation in \eqref{8psi8gl} generalizes
W.~N.~Bailey's~\cite[Eq.~(4.7)]{bail66} very-well-poised $_6\psi_6$ summation,
to which it reduces for $c\to 0$.

By an elementary computation it follows that H.~S.~Shukla's $_8\psi_8$
summation can also be written in the following eqivalent form:
\begin{multline}\label{88gl}
{}_8\psi_8\!\left[\begin{matrix}q\sqrt{a},-q\sqrt{a},b,c,aq^2/c,e,f,g\\
\sqrt{a},-\sqrt{a},aq/b,aq/c,c/q,aq/e,aq/f,aq/g\end{matrix}\,
;q,\frac{a^2}{befg}\right]\\
=\left(1-\frac{(1-a/bg)(1-cf/aq)(1-ce/aq)}{(1-c/gq)(1-ef/a)(1-c/bq)}\right)
\frac{(1-c/bq)(1-c/gq)}{(1-c/aq)(1-c/q)}\\\times
\frac{(q,aq,q/a,aq/be,aq/bf,aq/bg,a/ef,aq/eg,aq/fg;q)_{\infty}}
{(aq/b,aq/e,aq/f,aq/g,q/b,q/e,q/f,q/g,a^2/befg;q)_{\infty}},
\end{multline}
provided $|q|<1$ and $|a^2/befg|<1$.
In Section~\ref{sec8}, we provide a multidimensional generalization of this
form of H.~S.~Shukla's $_8\psi_8$ summation, see Theorem~\ref{m88}.

H.~S.~Shukla~\cite{shukla} derived \eqref{8psi8gl} by specializing a
transformation of a very-well-poised $_8\psi_8$ into a sum of three
balanced $_4\phi_3$ series due to M.~Jackson~\cite[Eq.~(3.1)]{mjackson8}.
Unlike H.~S.~Shukla~\cite{shukla}, we give a
proof of the $_8\psi_8$ summation formula \eqref{88gl} using
a weaker result, namely
L.~J.~Roger's~\cite{rogers} $_6\phi_5$ summation, which is a special
case of \eqref{88gl}. As further ingredients in our derivation of \eqref{88gl} 
we utilize a simple decomposition identity, see Eq.~\eqref{pbz2},
and an application of M.~E.~H.~Ismail's~\cite{ismail}
analytic continuation argument. We display our proof of \eqref{88gl} in
Section~\ref{sec1}.

Surprisingly, the whole analysis carries over to the multivariate case,
in the setting of {\em multiple basic hypergeometric series very-well-poised
over the root-system $A_{r-1}$}. In fact, the main achievement of this article
is an $A_{r-1}$ generalization of H.~S.~Shukla's $_8\psi_8$ summation
\eqref{88gl}, see Theorem~\ref{m88}.

Our article is organized as follows. In Section~\ref{sec1},
after briefly explaining some basic concepts which we need from
the theory of basic hypergeometric series~\cite{grhyp},
we give a new proof of H.~S.~Shukla's $_8\psi_8$ summation.
The one-dimensional analysis in Section~\ref{sec1} turns out to be
very much motivated by the multivariate case. In Section~\ref{sec8},
after some preparations, we state and prove our $A_{r-1}$ generalization
of H.~S.~Shukla's very-well-poised $_8\psi_8$ summation theorem.
Our $A_{r-1}$ $_8\psi_8$ summation in Theorem~\ref{m88} includes
R.~A.~Gustafson's~\cite{gusmult} $A_{r-1}$ $_6\psi_6$ summation as a
special case. Our proof utilizes an $A_{r-1}$ nonterminating
very-well-poised $_6\phi_5$ summation by S.~C.~Milne~\cite{milnun},
a partial fraction decomposition (see Lemma~\ref{pbz}),
and analytic continuation.

\section{Some basic concepts and proof of the $_8\psi_8$ summation
in one-dimension}\label{sec1}

\subsection{Some basic concepts}

We first recall some basic concepts from the theory of
basic hypergeometric series (cf.\ G.~Gasper and M.~Rahman~\cite{grhyp}).

By definition \eqref{qpochdef} it is clear that
$(q;q)_k^{-1}=0$, for $k=-1, -2, \ldots$. More generally,
$(q^{1+m};q)_k^{-1}=0$ for $k=-m-1,-m-2, \ldots$. Thus, using the
definitions \eqref{defhyp} and \eqref{defhypb},
a bilateral ${}_t\psi_t$ series becomes a
unilateral $_t\phi_{t-1}$ series if one of the lower parameters, say
$b_t$, is $q$ (or more generally, $q^j$ where $j$ is a positive integer).
In this case, the ${}_t\psi_t$ series terminates naturally from below.
Similarly, a $_t\phi_{t-1}$ series terminates naturally from above
if one of the upper parameters, say $a_t$, equals $q^{-n}$, $n=0,1,2,\ldots$.

The ratio test gives simple criteria of when the series in
\eqref{defhyp} and \eqref{defhypb} converge,
if they do not terminate. Remember that we assume $|q|<1$.
The $_t\phi_{t-1}$ series in \eqref{defhyp} converges absolutely in the radius
$|z|<1$, while the ${}_t\psi_t$ series in \eqref{defhypb} converges absolutely
in the annulus $|b_1\dots b_t/a_1\dots a_t|<|z|<1$. 

The classical theory of basic hypergeometric series consists of
several summation and transformation formulae involving $_t\phi_{t-1}$
or ${}_t\psi_t$ series. Some classical summation theorems for these
series require that the parameters satisfy the condition of being
very-well-poised. A $_t\phi_{t-1}$ basic hypergeometric series is called
{\em well-poised} if $a_1q=a_2b_1=\cdots=a_tb_{t-1}$.
It is called {\em very-well-poised} if it is well-poised and if
$a_2=q\sqrt{a_1}$ and $a_3=-q\sqrt{a_1}$.
Note that  the factor
\begin{equation}\label{vwp}
\frac{(q\sqrt{a_1},-q\sqrt{a_1};q)_k}{(\sqrt{a_1},-\sqrt{a_1};q)_k}=
\frac {1-a_1q^{2k}}{1-a_1}
\end{equation}
appears in a very-well-poised series.
The parameter $a_1$ is usually referred to as the
{\em special parameter} of such a series, and we call \eqref{vwp}
the {\em very-well-poised term} of the series.
Similarly, a bilateral $_t\psi_t$ basic
hypergeometric series is well-poised if
$a_1b_1=a_2b_2\cdots=a_tb_t$ and very-well-poised if, in addition,
$a_1=-a_2=qb_1=-qb_2$.

In our subsequent computations (in this section and in Section~\ref{sec8}),
we make heavily use of some elementary identities involving
$q$-shifted factorials, listed in
G.~Gasper and M.~Rahman~\cite[Appendix~I]{grhyp}.

\subsection{Proof of the
$\hbox{${}_{\boldsymbol 8}\boldsymbol\psi_{\boldsymbol 8}$}$ summation}

The main ingredient in our derivation of H.~S.~Shukla's $_8\psi_8$ summation
\eqref{88gl} is L.~J.~Roger's~\cite[p.~29, second eq.]{rogers}
nonterminating very-well-poised $_6\phi_5$ summation:
\begin{equation}\label{65gl}
{}_6\phi_5\!\left[\begin{matrix}a,\,q\sqrt{a},-q\sqrt{a},b,c,d\\
\sqrt{a},-\sqrt{a},aq/b,aq/c,aq/d\end{matrix}\,;q,
\frac{aq}{bcd}\right]
=\frac {(aq,aq/bc,aq/bd,aq/cd;q)_{\infty}}
{(aq/b,aq/c,aq/d,aq/bcd;q)_{\infty}},
\end{equation}
provided $|q|<1$ and $|aq/bcd|<1$. Note that \eqref{65gl} is equivalent to
the special case $b\to a$, $c\to 0$ of \eqref{88gl}.

We derive the $_8\psi_8$ summation \eqref{88gl} in two steps.
In the first step, we establish the $b\to a$ case of \eqref{88gl}
by using L.~J.~Roger's $_6\phi_5$ summation \eqref{65gl} twice.
I.e., we first establish the unilateral summation
\begin{multline}\label{87gl}
{}_8\phi_7\!\left[\begin{matrix}a,q\sqrt{a},-q\sqrt{a},c,aq^2/c,e,f,g\\
\sqrt{a},-\sqrt{a},aq/c,c/q,aq/e,aq/f,aq/g\end{matrix}\,
;q,\frac{a}{efg}\right]\\
=\left(1-\frac{(1-1/g)(1-cf/aq)(1-ce/aq)}{(1-c/gq)(1-ef/a)(1-c/aq)}\right)
\frac{(1-c/gq)}{(1-c/q)}\\\times
\frac{(aq,a/ef,aq/eg,aq/fg;q)_{\infty}}
{(aq/e,aq/f,aq/g,a/efg;q)_{\infty}},
\end{multline}
where $|q|<1$ and $|a/efg|<1$.
In the second step, we extend \eqref{87gl} to \eqref{88gl} by analytic
continuation, or equivalently, by an application of
M.~E.~H.~Ismail's~\cite{ismail} argument.

The details of the first step are as follows.
Since
\begin{equation}\label{pbz2}
\frac{(1-cq^{k-1})(1-aq^{k+1}/c)}{(1-c/q)(1-aq/c)}=
q^k+\frac{(1-aq^k)(1-q^k)}{(1-c/q)(1-aq/c)},
\end{equation}
we have
\begin{multline*}
{}_8\phi_7\!\left[\begin{matrix}a,q\sqrt{a},-q\sqrt{a},c,aq^2/c,e,f,g\\
\sqrt{a},-\sqrt{a},aq/c,c/q,aq/e,aq/f,aq/g\end{matrix}\,
;q,\frac{a}{efg}\right]\\
=\sum_{k=0}^{\infty}\frac{(1-aq^{2k})}{(1-a)}
\frac{(a,e,f,g;q)_k}{(q,aq/e,aq/f,aq/g;q)_k}\left(\frac a{efg}\right)^k
\frac{(1-cq^{k-1})(1-aq^{k+1}/c)}{(1-c/q)(1-aq/c)}\\
=\sum_{k=0}^{\infty}\frac{(1-aq^{2k})}{(1-a)}
\frac{(a,e,f,g;q)_k}{(q,aq/e,aq/f,aq/g;q)_k}\left(\frac{aq}{efg}\right)^k\\
+\sum_{k=0}^{\infty}\frac{(1-aq^{2k})}{(1-a)}
\frac{(a,e,f,g;q)_k}{(q,aq/e,aq/f,aq/g;q)_k}\left(\frac a{efg}\right)^k
\frac{(1-aq^k)(1-q^k)}{(1-c/q)(1-aq/c)}.
\end{multline*}
Now in the second sum, because of the factor $(1-q^k)$ in the numerator of the
summand, we shift the index $k\mapsto k+1$. We then obtain
\begin{multline*}
\sum_{k=0}^{\infty}\frac{(1-aq^{2k})}{(1-a)}
\frac{(a,e,f,g;q)_k}{(q,aq/e,aq/f,aq/g;q)_k}\left(\frac{aq}{efg}\right)^k\\
+\frac{a(1-aq)(1-aq^2)(1-e)(1-f)(1-g)}
{efg(1-c/q)(1-aq/c)(1-aq/e)(1-aq/f)(1-aq/g)}\\\times
\sum_{k=0}^{\infty}\frac{(1-aq^{2+2k})}{(1-aq^2)}
\frac{(aq^2,eq,fq,gq;q)_k}{(q,aq^2/e,aq^2/f,aq^2/g;q)_k}
\left(\frac a{efg}\right)^k.
\end{multline*}
Next, we simplify both sums by the $_6\phi_5$ summation in \eqref{65gl}
and obtain
\begin{multline}\label{gf1}
\frac{(aq,aq/ef,aq/eg,aq/fg;q)_{\infty}}{(aq/e,aq/f,aq/g,aq/efg;q)_{\infty}}\\
+\frac{a(1-aq)(1-aq^2)(1-e)(1-f)(1-g)}
{efg(1-c/q)(1-aq/c)(1-aq/e)(1-aq/f)(1-aq/g)}\\\times
\frac{(aq^3,aq/ef,aq/eg,aq/fg;q)_{\infty}}
{(aq^2/e,aq^2/f,aq^2/g,a/efg;q)_{\infty}}\\
=\left(1-\frac{(1-e)(1-f)(1-g)}{(1-c/q)(1-aq/c)(1-efg/a)}\right)
\frac{(aq,aq/ef,aq/eg,aq/fg;q)_{\infty}}{(aq/e,aq/f,aq/g,aq/efg;q)_{\infty}}.
\end{multline}
Since
\begin{multline*}
\left(1-\frac{(1-e)(1-f)(1-g)}{(1-c/q)(1-aq/c)(1-efg/a)}\right)
\frac{(1-a/efg)}{(1-a/ef)}\\
=\left(1-\frac{(1-1/g)(1-cf/aq)(1-ce/aq)}{(1-c/gq)(1-ef/a)(1-c/aq)}\right)
\frac{(1-c/gq)}{(1-c/q)}
\end{multline*}
(as can be readily checked by using a symbolic computer algebra program
such as {\em Maple} or {\em Mathematica}), the last expression in
\eqref{gf1} is equivalent to
\begin{equation*}
\left(1-\frac{(1-1/g)(1-cf/aq)(1-ce/aq)}{(1-c/gq)(1-ef/a)(1-c/aq)}\right)
\frac{(1-c/gq)}{(1-c/q)}
\frac{(aq,a/ef,aq/eg,aq/fg;q)_{\infty}}
{(aq/e,aq/f,aq/g,a/efg;q)_{\infty}},
\end{equation*}
which is the right side of \eqref{87gl}.

Having established \eqref{87gl}, we are now ready to proceed
with the second step, where we extend the unilateral summation \eqref{87gl} to
the bilateral \eqref{88gl} by
analytic continuation,  by a method commonly referred to as
``Ismail's argument" (see M.~E.~H.~Ismail~\cite{ismail}, and
R.~Askey and M.~E.~H.~Ismail~\cite{askmail}). This works as follows:
Both sides of the identity in \eqref{88gl} are analytic in $b^{-1}$ in a domain
around the origin. Now, the identity is true for $b=aq^{-m}$, for all
$m=0,1,2,\ldots$, by the $_8\phi_7$ summation in \eqref{87gl}
(see below for the details). Since $\lim_{m\to\infty}q^m/a=0$ is an interior
point in the domain of analyticity of $b^{-1}$, by the identity
theorem of analytic functions,
we establish the identity \eqref{88gl} for $b^{-1}$ throughout the whole
domain. Finally, by analytic continuation we esatblish the identity
\eqref{88gl} to be valid for $|b^{-1}|<|efg/a^2|$, the region of convergence
of the series.

We still need to show that the identity \eqref{88gl} is true
when $b=aq^{-m}$. In this case, the left side
of \eqref{88gl} is
\begin{equation}\label{longl1}
\sum_{k=-m}^{\infty}\frac{(1-aq^{2k})}{(1-a)}
\frac{(aq^{-m},c,aq^2/c,e,f,g;q)_k}{(q^{1+m},aq/c,c/q,aq/e,aq/f,aq/g;q)_k}
\left(\frac{aq^m}{efg}\right)^k.
\end{equation}
We shift the summation index in \eqref{longl1} by
$k\mapsto k-m$ and obtain
\begin{multline*}
\frac{(1-aq^{-2m})}{(1-a)}
\frac{(aq^{-m},c,aq^2/c,e,f,g;q)_{-m}}
{(q^{1+m},aq/c,c/q,aq/e,aq/f,aq/g;q)_{-m}}
\left(\frac{aq^m}{efg}\right)^{-m}\\\times
\sum_{k=0}^{\infty}
\frac{(1-aq^{-2m+2k})(aq^{-2m},cq^{-m},aq^{2-m}/c,eq^{-m},fq^{-m},gq^{-m};q)_k}
{(1-a^{-2m})(q,aq^{1-m}/c,cq^{-m-1},aq^{1-m}/e,aq^{1-m}/f,aq^{1-m}/g;q)_k}
\!\left(\frac{aq^m}{efg}\right)^k\!\!.
\end{multline*}
Now we apply the $a\mapsto aq^{-2m}$, $c\mapsto cq^{-m}$,
$e\mapsto eq^{-m}$, $f\mapsto fq^{-m}$, and $g\mapsto gq^{-m}$ case of
the summation formula in \eqref{87gl} to simplify this expression to
\begin{multline*}
\frac{(1-aq^{-2m})}{(1-a)}
\frac{(aq^{-m},c,aq^2/c,e,f,g;q)_{-m}}
{(q^{1+m},aq/c,c/q,aq/e,aq/f,aq/g;q)_{-m}}
\left(\frac{aq^m}{efg}\right)^{-m}\\\times
\left(1-\frac{(1-q^m/g)(1-cf/aq)(1-ce/aq)}
{(1-c/gq)(1-ef/a)(1-cq^{m}/aq)}\right)
\frac{(1-c/gq)}{(1-cq^{-m-1})}\\\times
\frac{(aq^{1-2m},a/ef,aq/eg,aq/fg;q)_{\infty}}
{(aq^{1-m}/e,aq^{1-m}/f,aq^{1-m}/g,aq^m/efg;q)_{\infty}}.
\end{multline*}
Now, this can easily be further transformed into
\begin{multline*}
\left(1-\frac{(1-q^m/g)(1-cf/aq)(1-ce/aq)}{(1-c/gq)(1-ef/a)(1-cq^m/aq)}\right)
\frac{(1-cq^m/aq)(1-c/gq)}{(1-c/aq)(1-c/q)}\\\times
\frac{(q,aq,q/a,q^{1+m}/e,q^{1+m}/f,q^{1+m}/g,a/ef,aq/eg,aq/fg;q)_{\infty}}
{(q^{1+m},aq/e,aq/f,aq/g,q^{1+m}/a,q/e,q/f,q/g,aq^{m}/efg;q)_{\infty}},
\end{multline*}
which is exactly the $b=aq^{-m}$ case of the right side of \eqref{88gl}.
\qed

\section{An $A_{r-1}$ very-well-poised $_8\psi_8$ summation formula}\label{sec8}

\subsection{Preliminaries on
$\hbox{$\boldsymbol A_{{\boldsymbol r}{\boldsymbol -}{\boldsymbol 1}}$}$
basic hypergeometric series}

We consider multiple series of the form
\begin{equation}\label{skgl}
\sum_{k_1,\dots,k_r=-\infty}^{\infty}S({\mathbf k}),
\end{equation}
where ${\mathbf k}=(k_1,\dots,k_r)$, which reduce to classical
(bilateral) basic hypergeometric series when $r=1$.
We call such a multiple basic hypergeometric series {\em well-poised}
if it reduces to a well-poised series when $r=1$. {\em Very-well-poised}
multiple basic hypergeometric series are defined analogously.
In case these series do not terminate from below, we also call such series
{\em multilateral} basic hypergeometric series.

In our particular cases, we also have
\begin{equation}\label{arvandy}
\prod_{1\le i<j<r} \left(\frac {z_iq^{k_i}-z_jq^{k_j}} {z_i-z_j}\right)
\end{equation} 
as a factor of $S({\mathbf k})$. A typical example
is the right side of \eqref{r66def}. Since we may associate \eqref{arvandy}
with the product side of the Weyl denominator formula for the
root system $A_{r-1}$ (see e.g.\ D.~Stanton~\cite{stan}), we call our series
$A_{r-1}$ basic hypergeometric series, in accordance with I.~M.~Gessel and
C.~Krattenthaler~\cite[Eq.~(7.1)]{geskratt}.
Very often these series are also called $U\!(r)$ basic hypergeometric series,
where $U\!(r)$ is the unitary group.
For some selected results in the theory of $A_{r-1}$ basic hypergeometric
series, see the references \cite{bhatmil,bhatms,geskratt,gusmult,gus,
milqana,milnun,Milgauss,milne,milschloss,schlossmmi,schlnammi,schlidp}.

For convenience, we frequently use the notation
$|{\mathbf k}|:= k_1+\dots+k_r$. Furthermore,
we often use capital letters to abbreviate the ($r$-fold)
products of certain variables. Specifically, in this article we use
$B\equiv b_1\cdots b_r$, $E\equiv e_1\cdots e_r$, and $Y\equiv y_1\cdots y_r$,
respectively. 

Since multidimensional $_6\phi_5$ series play a significant role in our
derivation of our $A_{r-1}$ $_8\psi_8$ summation, we find it useful to make
the following definition.
Let $a$, $b_1,\dots,b_r$, $c$, $d$, $z_1,\dots,z_r$,
and $w$ be indeterminate.
We define for $r\ge 1$,
\begin{multline}\label{r66def}
{}_6\Phi_5^{(r)}\!\left[a;b_1,\dots,b_r;c,d;z_1,\dots,z_r
\big|\,q,w\right]\\
:=\sum_{k_1,\dots,k_r=0}^{\infty}
\Bigg(\prod_{1\le i<j\le r}
\left(\frac{z_iq^{k_i}-z_jq^{k_j}}{z_i-z_j}\right)
\prod_{i=1}^r\left(\frac{1-az_iq^{k_i+|{\mathbf k}|}}{1-az_i}\right)
\prod_{i,j=1}^r\frac{(b_jz_i/z_j;q)_{k_i}}{(qz_i/z_j;q)_{k_i}}\\\times
\prod_{i=1}^r\frac{(az_i;q)_{|{\mathbf k}|}\,(cz_i;q)_{k_i}}
{(az_iq/b_i;q)_{|{\mathbf k}|}\,(az_iq/d;q)_{k_i}}\cdot
\frac{(d;q)_{|{\mathbf k}|}}{(aq/c;q)_{|{\mathbf k}|}}\,
w^{|{\mathbf k}|}\Bigg).
\end{multline}

The above $_6\Phi_5^{(r)}$ series is an $r$-dimensional $_6\phi_5$ series
(which reduces to a classical very-well-posied $_6\phi_5$ when $r=1$).

In our proof of Theorem~\ref{m88}, or more precisely, of the intermediate
Proposition~\ref{p88}, we utilize
S.~C.~Milne's~\cite[Theorem~1.44]{milnun} $A_{r-1}$ extension of
L.~J.~Roger's $_6\phi_5$ summation theorem.

\begin{Theorem}[(Milne) An $A_{r-1}$ nonterminating very-well-poised
$_6\phi_5$ summation]\label{r65}
Let $a$, $b_1,\dots,b_r$, $c$, $d$, and $z_1,\dots,z_r$, be indeterminate,
let $B\equiv b_1\cdots b_r$, $r\ge 1$, and suppose that
none of the denominators in \eqref{r65gl} vanishes. Then
\begin{multline}\label{r65gl}
{}_6\Phi_5^{(r)}\!\left[a;b_1,\dots,b_r;c,d;z_1,\dots,z_r
\big|\,q,\frac{aq}{Bcd}\right]\\
=\frac{(aq/Bc,aq/cd;q)_{\infty}}{(aq/Bcd,aq/c;q)_{\infty}}
\prod_{i=1}^r\frac{(az_iq,az_iq/b_id;q)_{\infty}}
{(az_iq/d,az_iq/b_i;q)_{\infty}},
\end{multline}
provided $|q|<1$ and $|aq/Bcd|<1$.
\end{Theorem}
The $r=1$ case of \eqref{r65gl} clearly reduces to \eqref{65gl}.

Further, we make use of the following ($q$-analogue of the)
{\em partial fraction decomposition}
\begin{equation}\label{pfd}
\prod_{i=1}^r\frac {(1-tz_iy_i)} {(1-tz_i)}=y_1y_2\dots y_r+
\sum_{l=1}^r\frac {\prod_{i=1}^r(1-y_iz_i/z_l)}
{(1-tz_l)\prod_{\begin{smallmatrix}i=1\\i\neq l\end{smallmatrix}}^r
(1-z_i/z_l)}
\end{equation}
(see~\cite[Appendix]{Milgauss}). In particular, for a multivariate extension
of \eqref{pbz2}, we utilize the following extension of \eqref{pfd}:
\begin{Lemma}\label{pbz}
Let $Y\equiv y_1y_2\dots y_r$. Then
\begin{equation*}
\frac{(1-uY)}{(1-u)}
\prod_{i=1}^r\frac {(1-tz_iy_i)} {(1-tz_i)}=Y+
\sum_{l=1}^r\frac {(1-uYtz_l)\prod_{i=1}^r(1-y_iz_i/z_l)}
{(1-u)(1-tz_l)\prod_{\begin{smallmatrix}i=1\\i\neq l\end{smallmatrix}}^r
(1-z_i/z_l)}.
\end{equation*}
\end{Lemma}
\begin{proof}
Since
\begin{equation*}
1-Y=\sum_{l=1}^r\frac {\prod_{i=1}^r(1-y_iz_i/z_l)}
{\prod_{\begin{smallmatrix}i=1\\i\neq l\end{smallmatrix}}^r
(1-z_i/z_l)},
\end{equation*}
by the $t=0$ case of \eqref{pfd}, we have
\begin{multline*}
\frac{(1-uY)}{(1-u)}
\prod_{i=1}^r\frac {(1-tz_iy_i)} {(1-tz_i)}\\
=\left(1+\frac{(1-Y)}{(1-u)}u\right)
\left(Y+\sum_{l=1}^r\frac {\prod_{i=1}^r(1-y_iz_i/z_l)}
{(1-tz_l)\prod_{\begin{smallmatrix}i=1\\i\neq l\end{smallmatrix}}^r
(1-z_i/z_l)}\right)\\
=Y+\sum_{l=1}^r\left(\frac{(1-u+(1-tz_l)uY+(1-Y)u)}{(1-u)}\cdot
\frac {\prod_{i=1}^r(1-y_iz_i/z_l)}
{(1-tz_l)\prod_{\begin{smallmatrix}i=1\\i\neq l\end{smallmatrix}}^r
(1-z_i/z_l)}\right)\\
=Y+\sum_{l=1}^r\frac {(1-uYtz_l)\prod_{i=1}^r(1-y_iz_i/z_l)}
{(1-u)(1-tz_l)\prod_{\begin{smallmatrix}i=1\\i\neq l\end{smallmatrix}}^r
(1-z_i/z_l)}.
\end{multline*}
\end{proof}

\begin{Remark}
In the multivariate analysis of our proof of Theorem~\ref{m88},
the partial fraction decomposition of Lemma~\ref{pbz} plays a crucial role.
Applications of partial fraction decomposition have often proved to be
useful in the derivation of results for $A_{r-1}$ series, see e.g.
\cite{bhatmil,gus,milqana,Milgauss,schlossmmi,schlnammi}.
\end{Remark}

After these preparations, we are ready to state and prove our multiple
extension of \eqref{88gl}.

\subsection{The main result}

Our $r$-dimensional generalization of
H.~S.~Shukla's~\cite[Eq.~(4.1)]{shukla} very-well-poised
$_8\psi_8$ summation is as follows:

\begin{Theorem}[An $A_{r-1}$ very-well-poised $_8\psi_8$ summation]\label{m88}
Let $a$, $b_1,\dots,b_r$, $c$, $e_1,\dots,e_r$, $f$, $g$, and $z_1,\dots,z_r$
be indeterminate, let $B\equiv b_1\cdots b_r$, $E\equiv e_1\cdots e_r$,
$r\ge 1$, and
suppose that none of the denominators in \eqref{m88gl} vanishes. Then
\begin{multline}\label{m88gl}
\sum_{k_1,\dots,k_r=-\infty}^{\infty}
\Bigg(\prod_{1\le i<j\le r}
\left(\frac{z_iq^{k_i}-z_jq^{k_j}}{z_i-z_j}\right)
\prod_{i=1}^r\left(\frac{1-az_iq^{k_i+|{\mathbf k}|}}{1-az_i}\right)\\\times
\prod_{i,j=1}^r\frac{(e_jz_i/z_j;q)_{k_i}}{(az_iq/b_jz_j;q)_{k_i}}
\prod_{i=1}^r\frac{(b_iz_i;q)_{|{\mathbf k}|}\,(fz_i;q)_{k_i}}
{(az_iq/e_i;q)_{|{\mathbf k}|}\,(az_iq/g;q)_{k_i}}\cdot
\frac{(g;q)_{|{\mathbf k}|}}{(aq/f;q)_{|{\mathbf k}|}}\\\times
\frac{(1-cq^{|{\mathbf k}|-1})}{(1-c/q)}
\prod_{i=1}^r\frac{(1-az_iq^{k_i+1}/c)}{(1-az_iq/c)}\cdot
\left(\frac{a^{r+1}}{BEfg}\right)^{|{\mathbf k}|}\Bigg)\\
=\left(1-\frac{(1-a^r/Bg)(1-cf/aq)}{(1-c/gq)(1-Ef/a)}
\prod_{i=1}^r\frac{(1-ce_i/az_iq)}{(1-c/b_iz_iq)}\right)
\prod_{i=1}^r\frac{(1-c/b_iz_iq)}{(1-c/az_iq)}\\\times
\frac{(1-c/gq)}{(1-c/q)}
\frac{(a/Ef,aq/fg,a^rq/Bg;q)_{\infty}}{(a^{r+1}/BEfg,aq/f,q/g;q)_{\infty}}
\prod_{i,j=1}^r\frac{(qz_i/z_j,az_iq/b_je_iz_j;q)_{\infty}}
{(az_iq/b_jz_j,z_iq/e_iz_j;q)_{\infty}}\\\times
\prod_{i=1}^r\frac{(az_iq,q/az_i,az_iq/e_ig,aq/b_ifz_i;q)_{\infty}}
{(az_iq/e_i,az_iq/g,q/b_iz_i,q/fz_i;q)_{\infty}},
\end{multline}
provided $|q|<1$ and $|a^{r+1}/BEfg|<1$.
\end{Theorem}

Theorem~\ref{m88} generalizes R.~A.~Gustafson's~\cite[Theorem~1.15]{gusmult}
$A_{r-1}$ $_6\psi_6$ summation, to which it reduces for $c\to 0$.

Following closely the univariate analysis of Section~\ref{sec1},
we prove Theorem~\ref{m88} in two steps. First, we prove the
$b_i=a$, $i=1,\dots,r$, special case of Theorem~\ref{m88}, which is
Proposition~\ref{p88} below. Then we extend Proposition~\ref{p88} to
Theorem~\ref{m88} by an $r$-fold application of
M.~E.~H.~Ismail's~\cite{ismail} analytic continuation argument.

\begin{Proposition}[An $A_{r-1}$ nonterminating very-well-poised $_8\phi_7$
summation]\label{p88}
Let $a$, $c$, $e_1,\dots,e_r$, $f$, $g$, and $z_1,\dots,z_r$
be indeterminate, let $E\equiv e_1\cdots e_r$, $r\ge 1$, and
suppose that none of the denominators in \eqref{p88gl} vanishes. Then
\begin{multline}\label{p88gl}
\sum_{k_1,\dots,k_r=0}^{\infty}
\Bigg(\prod_{1\le i<j\le r}
\left(\frac{z_iq^{k_i}-z_jq^{k_j}}{z_i-z_j}\right)
\prod_{i=1}^r\left(\frac{1-az_iq^{k_i+|{\mathbf k}|}}{1-az_i}\right)\\\times
\prod_{i,j=1}^r\frac{(e_jz_i/z_j;q)_{k_i}}{(qz_i/z_j;q)_{k_i}}
\prod_{i=1}^r\frac{(az_i;q)_{|{\mathbf k}|}\,(fz_i;q)_{k_i}}
{(az_iq/e_i;q)_{|{\mathbf k}|}\,(az_iq/g;q)_{k_i}}\cdot
\frac{(g;q)_{|{\mathbf k}|}}{(aq/f;q)_{|{\mathbf k}|}}\\\times
\frac{(1-cq^{|{\mathbf k}|-1})}{(1-c/q)}
\prod_{i=1}^r\frac{(1-az_iq^{k_i+1}/c)}{(1-az_iq/c)}\cdot
\left(\frac{a}{Efg}\right)^{|{\mathbf k}|}\Bigg)\\
=\left(1-\frac{(1-1/g)(1-cf/aq)}{(1-c/gq)(1-Ef/a)}
\prod_{i=1}^r\frac{(1-ce_i/az_iq)}{(1-c/az_iq)}\right)
\frac{(1-c/gq)}{(1-c/q)}\\\times
\frac{(a/Ef,aq/fg;q)_{\infty}}{(a/Efg,aq/f;q)_{\infty}}
\prod_{i=1}^r\frac{(az_iq,az_iq/e_ig;q)_{\infty}}
{(az_iq/e_i,az_iq/g;q)_{\infty}},
\end{multline}
provided $|q|<1$ and $|a/Efg|<1$.
\end{Proposition}
Proposition~\ref{p88} generalizes S.~C.~Milne's $A_{r-1}$ $_6\phi_5$
summation in Theorem~\ref{r65}, to which it reduces for $c=0$.
\begin{proof}[Proof of Proposition~\ref{p88}]
Since
\begin{multline*}
\frac{(1-cq^{|{\mathbf k}|-1})}{(1-c/q)}
\prod_{i=1}^r\frac{(1-az_iq^{k_i+1}/c)}{(1-az_iq/c)}\\=
q^{|{\mathbf k}|}+
\sum_{l=1}^r\frac {(1-az_lq^{|{\mathbf k}|})\prod_{i=1}^r(1-q^{k_i}z_i/z_l)}
{(1-c/q)(1-az_lq/c)\prod_{\begin{smallmatrix}i=1\\i\neq l\end{smallmatrix}}^r
(1-z_i/z_l)},
\end{multline*}
by the $t\mapsto aq/c$, $u\mapsto c/q$, $y_i\mapsto q^{k_i}$,
$i=1,\dots,r$, case of Lemma~\ref{pbz}, we have
\begin{multline*}
\sum_{k_1,\dots,k_r=0}^{\infty}
\Bigg(\prod_{1\le i<j\le r}
\left(\frac{z_iq^{k_i}-z_jq^{k_j}}{z_i-z_j}\right)
\prod_{i=1}^r\left(\frac{1-az_iq^{k_i+|{\mathbf k}|}}{1-az_i}\right)\\\times
\prod_{i,j=1}^r\frac{(e_jz_i/z_j;q)_{k_i}}{(qz_i/z_j;q)_{k_i}}
\prod_{i=1}^r\frac{(az_i;q)_{|{\mathbf k}|}\,(fz_i;q)_{k_i}}
{(az_iq/e_i;q)_{|{\mathbf k}|}\,(az_iq/g;q)_{k_i}}\cdot
\frac{(g;q)_{|{\mathbf k}|}}{(aq/f;q)_{|{\mathbf k}|}}\\\times
\left(\frac{a}{Efg}\right)^{|{\mathbf k}|}
\frac{(1-cq^{|{\mathbf k}|-1})}{(1-c/q)}
\prod_{i=1}^r\frac{(1-az_iq^{k_i+1}/c)}{(1-az_iq/c)}\Bigg)\\
=\sum_{k_1,\dots,k_r=0}^{\infty}\Bigg(\prod_{1\le i<j\le r}
\left(\frac{z_iq^{k_i}-z_jq^{k_j}}{z_i-z_j}\right)
\prod_{i=1}^r\left(\frac{1-az_iq^{k_i+|{\mathbf k}|}}{1-az_i}\right)
\prod_{i,j=1}^r\frac{(e_jz_i/z_j;q)_{k_i}}{(qz_i/z_j;q)_{k_i}}\\\times
\prod_{i=1}^r\frac{(az_i;q)_{|{\mathbf k}|}\,(fz_i;q)_{k_i}}
{(az_iq/e_i;q)_{|{\mathbf k}|}\,(az_iq/g;q)_{k_i}}\cdot
\frac{(g;q)_{|{\mathbf k}|}}{(aq/f;q)_{|{\mathbf k}|}}
\left(\frac{aq}{Efg}\right)^{|{\mathbf k}|}\Bigg)\\
+\sum_{l=1}^r
\sum_{k_1,\dots,k_r=0}^{\infty}
\Bigg(\prod_{1\le i<j\le r}
\left(\frac{z_iq^{k_i}-z_jq^{k_j}}{z_i-z_j}\right)
\prod_{i=1}^r\left(\frac{1-az_iq^{k_i+|{\mathbf k}|}}{1-az_i}\right)\\\times
\prod_{i,j=1}^r\frac{(e_jz_i/z_j;q)_{k_i}}{(qz_i/z_j;q)_{k_i}}
\prod_{i=1}^r\frac{(az_i;q)_{|{\mathbf k}|}\,(fz_i;q)_{k_i}}
{(az_iq/e_i;q)_{|{\mathbf k}|}\,(az_iq/g;q)_{k_i}}\cdot
\frac{(g;q)_{|{\mathbf k}|}}{(aq/f;q)_{|{\mathbf k}|}}\\\times
\left(\frac{a}{Efg}\right)^{|{\mathbf k}|}
\frac {(1-az_lq^{|{\mathbf k}|})\prod_{i=1}^r(1-q^{k_i}z_i/z_l)}
{(1-c/q)(1-az_lq/c)\prod_{\begin{smallmatrix}i=1\\i\neq l\end{smallmatrix}}^r
(1-z_i/z_l)}\Bigg).
\end{multline*}
We arrived at a sum of $1+r$ infinite multisums.
In the $(1+l)$th sum, for $l=1,\dots,r$, due to the factor
$\prod_{i=1}^r(1-q^{k_i}z_i/z_l)$ in the numerator of the summand, 
we shift the index $k_l\mapsto k_l+1$. We then obtain
\begin{multline*}
{}_6\Phi_5^{(r)}\!\left[a;e_1,\dots,e_r;f,g;z_1,\dots,z_r
\big|\,q,\frac{aq}{Efg}\right]\\
+\sum_{l=1}^r
\frac{a(1-fz_l)(1-g)(1-az_lq^2)}
{Efg(1-c/q)(1-az_lq/c)(1-az_lq/g)(1-aq/f)}\\\times
\frac{\prod_{i=1}^r(1-az_iq)\prod_{i=1}^r(1-e_iz_l/z_i)}
{\prod_{i=1}^r(1-az_iq/e_i)
\prod_{\begin{smallmatrix}i=1\\i\neq l\end{smallmatrix}}^r(1-z_l/z_i)}\\\times
{}_6\Phi_5^{(r)}\!\bigg[aq;e_1,\dots,e_{l-1},e_lq,e_{l+1},\dots,e_r;f,gq;\\
z_1,\dots,z_{l-1},z_lq,z_{l+1},\dots,z_r
\big|\,q,\frac{a}{Efg}\bigg].
\end{multline*}
Next, we simplify all the $1+r$ ${}_6\Phi_5^{(r)}$ series according to
Theorem~\ref{r65}, and obtain
\begin{multline}\label{p881gl}
\frac{(aq/Ef,aq/fg;q)_{\infty}}{(aq/Efg,aq/f;q)_{\infty}}
\prod_{i=1}^r\frac{(az_iq,az_iq/e_ig;q)_{\infty}}
{(az_iq/g,az_iq/e_i;q)_{\infty}}\\
+\sum_{l=1}^r
\frac{a(1-fz_l)(1-g)(1-az_lq^2)}
{Efg(1-c/q)(1-az_lq/c)(1-az_lq/g)(1-aq/f)}\\\times
\frac{\prod_{i=1}^r(1-az_iq)\prod_{i=1}^r(1-e_iz_l/z_i)}
{\prod_{i=1}^r(1-az_iq/e_i)
\prod_{\begin{smallmatrix}i=1\\i\neq l\end{smallmatrix}}^r(1-z_l/z_i)}\\\times
\frac{(aq/Ef,aq/fg;q)_{\infty}}{(a/Efg,aq^2/f;q)_{\infty}}
\frac{(1-az_lq/g)}{(1-az_lq^2)}
\prod_{i=1}^r\frac{(az_iq^2,az_iq/e_ig;q)_{\infty}}
{(az_iq/g,az_iq^2/e_i;q)_{\infty}}\\
=\left(1-\sum_{l=1}^r\frac{(1-fz_l)(1-g)\prod_{i=1}^r(1-e_iz_l/z_i)}
{(1-c/q)(1-az_lq/c)(1-Efg/a)
\prod_{\begin{smallmatrix}i=1\\i\neq l\end{smallmatrix}}^r(1-z_l/z_i)}
\right)\\\times
\frac{(aq/Ef,aq/fg;q)_{\infty}}{(aq/Efg,aq/f;q)_{\infty}}
\prod_{i=1}^r\frac{(az_iq,az_iq/e_ig;q)_{\infty}}
{(az_iq/g,az_iq/e_i;q)_{\infty}}.
\end{multline}
Now we apply the $t\mapsto c/aq$, $u\mapsto aq/cEf$, $y_i\mapsto e_i$,
and $z_i\mapsto 1/z_i$, $i=1,\dots,r$, case of Lemma~\ref{pbz},
which can be rewritten as
\begin{multline*}
\sum_{l=1}^r\frac {(1-fz_l)\prod_{i=1}^r(1-e_iz_l/z_i)}
{(1-cEf/aq)(1-az_lq/c)
\prod_{\begin{smallmatrix}i=1\\i\neq l\end{smallmatrix}}^r(1-z_l/z_i)}\\
=1-\frac{(1-cf/aq)}{(1-cEf/aq)}
\prod_{i=1}^r\frac {(1-ce_i/az_iq)} {(1-c/az_iq)},
\end{multline*}
to simplify the expression obtained in \eqref{p881gl} to
\begin{multline}\label{p882gl}
\left(1-\frac{(1-g)(1-cEf/aq)}{(1-c/q)(1-Efg/a)}\left(
1-\frac{(1-cf/aq)}{(1-cEf/aq)}
\prod_{i=1}^r\frac {(1-ce_i/az_iq)} {(1-c/az_iq)}\right)\right)\\\times
\frac{(aq/Ef,aq/fg;q)_{\infty}}{(aq/Efg,aq/f;q)_{\infty}}
\prod_{i=1}^r\frac{(az_iq,az_iq/e_ig;q)_{\infty}}
{(az_iq/g,az_iq/e_i;q)_{\infty}}.
\end{multline}
Finally, using
\begin{equation*}
1-\frac{(1-g)(1-cEf/aq)}{(1-c/q)(1-Efg/a)}=
\frac{(1-c/gq)(1-a/Ef)}{(1-c/q)(1-a/Efg)},
\end{equation*}
we can easily transform the expression in \eqref{p882gl} into
\begin{multline*}
\left(1-\frac{(1-1/g)(1-cf/aq)}{(1-c/gq)(1-Ef/a)}
\prod_{i=1}^r\frac{(1-ce_i/az_iq)}{(1-c/az_iq)}\right)
\frac{(1-c/gq)}{(1-c/q)}\\\times
\frac{(a/Ef,aq/fg;q)_{\infty}}{(a/Efg,aq/f;q)_{\infty}}
\prod_{i=1}^r\frac{(az_iq,az_iq/e_ig;q)_{\infty}}
{(az_iq/e_i,az_iq/g;q)_{\infty}},
\end{multline*}
which is the right side of \eqref{p88gl}, as desired.
\end{proof}

Similar to the one-dimensional case, where we deduced the bilateral summation
\eqref{88gl} from the unilateral summation \eqref{87gl} by using
M.~E.~H.~Ismail's~\cite{ismail} argument, we can now readily
deduce Theorem~\ref{m88} from Proposition~\ref{p88}.

\begin{proof}[Proof of Theorem~\ref{m88}]
To establish \eqref{m88gl}, we apply Ismail's argument successively to the
parameters $b_1^{-1},\dots,b_r^{-1}$ using Proposition~\ref{p88}.
Both sides of the multiple series identity in \eqref{m88gl} are
analytic in each of the parameters $b_1^{-1},\dots,b_r^{-1}$ in a domain
around the origin. Now, the identity is true for $b_1=aq^{-m_1},
b_2=aq^{-m_2},\dots,$ and $b_r=aq^{-m_r}$, by the $A_{r-1}$ summation
in Proposition~\ref{p88} (see below for the details). This holds for
all $m_1,\dots,m_r\ge 0$.
Since $\lim_{m_1\to\infty}q^{m_1}/a=0$ is an interior point
in the domain of analyticity of $b_1^{-1}$, by the identity theorem
of analytic functions, we obtain an identity for $b_1^{-1}$.
By iterating this argument for $b_2^{-1},\dots,b_r^{-1}$, and
analytic continuation, we establish
\eqref{m88gl} for general $b_1^{-1},\dots,b_r^{-1}$
where $|B^{-1}|<|Efg/a^{r+1}|$.

The details are displayed as follows. Setting $b_i=aq^{-m_i}$, for
$i=1,\dots,r$, the left side of \eqref{m88gl} becomes
\begin{multline}\label{longlm1}
\sum_{\begin{smallmatrix}-m_i\le k_i\le\infty\\
i=1,\dots,r\end{smallmatrix}}\Bigg(
\prod_{1\le i<j\le r}
\left(\frac{z_iq^{k_i}-z_jq^{k_j}}{z_i-z_j}\right)
\prod_{i=1}^r\left(\frac{1-az_iq^{k_i+|{\mathbf k}|}}{1-az_i}\right)\\\times
\prod_{i,j=1}^r\frac{(e_jz_i/z_j;q)_{k_i}}{(q^{1+m_j}z_i/z_j;q)_{k_i}}
\prod_{i=1}^r\frac{(az_iq^{-m_i};q)_{|{\mathbf k}|}\,(fz_i;q)_{k_i}}
{(az_iq/e_i;q)_{|{\mathbf k}|}\,(az_iq/g;q)_{k_i}}\\\times
\frac{(g;q)_{|{\mathbf k}|}}{(aq/f;q)_{|{\mathbf k}|}}
\frac{(1-cq^{|{\mathbf k}|-1})}{(1-c/q)}
\prod_{i=1}^r\frac{(1-az_iq^{k_i+1}/c)}{(1-az_iq/c)}\cdot
\left(\frac{aq^{|{\mathbf m}|}}{Efg}\right)^{|{\mathbf k}|}\Bigg).
\end{multline}
We shift the summation indices in \eqref{longlm1} by
$k_i\mapsto k_i-m_i$, for $i=1,\dots,r$, and obtain
\begin{multline}\label{longlm}
\prod_{1\le i<j\le r}
\left(\frac{z_iq^{-m_i}-z_jq^{-m_j}}{z_i-z_j}\right)
\prod_{i=1}^r\left(\frac{1-az_iq^{-m_i-|{\mathbf m}|}}{1-az_i}\right)\\\times
\prod_{i,j=1}^r\frac{(e_jz_i/z_j;q)_{-m_i}}{(q^{1+m_j}z_i/z_j;q)_{-m_i}}
\prod_{i=1}^r\frac{(az_iq^{-m_i};q)_{-|{\mathbf m}|}\,(fz_i;q)_{-m_i}}
{(az_iq/e_i;q)_{-|{\mathbf m}|}\,(az_iq/g;q)_{-m_i}}\\\times
\frac{(g;q)_{-|{\mathbf m}|}}{(aq/f;q)_{-|{\mathbf m}|}}
\frac{(1-cq^{-|{\mathbf m}|-1})}{(1-c/q)}
\prod_{i=1}^r\frac{(1-az_iq^{1-m_i}/c)}{(1-az_iq/c)}\cdot
\left(\frac{aq^{|{\mathbf m}|}}{Efg}\right)^{-|{\mathbf m}|}\\\times
\sum_{\begin{smallmatrix}-m_i\le k_i\le\infty\\
i=1,\dots,r\end{smallmatrix}}\Bigg(
\prod_{1\le i<j\le r}
\left(\frac{z_iq^{-m_i+k_i}-z_jq^{-m_j+k_j}}{z_iq^{-m_i}-z_jq^{-m_j}}\right)
\prod_{i=1}^r\left(\frac{1-az_iq^{-m_i-|{\mathbf m}|+k_i+|{\mathbf k}|}}
{1-az_iq^{-m_i-|{\mathbf m}|}}\right)\\\times
\prod_{i,j=1}^r\frac{(e_jq^{-m_i}z_i/z_j;q)_{k_i}}
{(q^{1+m_j-m_i}z_i/z_j;q)_{k_i}}
\prod_{i=1}^r\frac{(az_iq^{-m_i-|{\mathbf m}|};q)_{|{\mathbf k}|}\,
(fz_iq^{-m_i};q)_{k_i}}
{(az_iq^{1-|{\mathbf m}|}/e_i;q)_{|{\mathbf k}|}\,(az_iq^{1-m_i}/g;q)_{k_i}}
\\\times\frac{(gq^{-|{\mathbf m}|};q)_{|{\mathbf k}|}}
{(aq^{1-|{\mathbf m}|}/f;q)_{|{\mathbf k}|}}
\frac{(1-cq^{|{\mathbf k}|-|{\mathbf m}|-1})}{(1-cq^{-|{\mathbf m}|-1})}
\prod_{i=1}^r\frac{(1-az_iq^{k_i+1-m_i}/c)}{(1-az_iq^{1-m_i}/c)}\cdot
\left(\frac{aq^{|{\mathbf m}|}}{Efg}\right)^{|{\mathbf k}|}\Bigg).
\end{multline}
Next, using the identities
\begin{equation*}
\prod_{i,j=1}^r(e_jz_i/z_j;q)_{-m_i}=
(-1)^{r|\mathbf m|}E^{-|\mathbf m|}
q^{r\sum_{i=1} ^r\binom{m_i+1}2}\prod_{i=1}^rz_i^{|{\mathbf m}|-rm_i}
\prod_{i,j=1}^r(z_iq/e_iz_j;q)_{m_j}^{-1},
\end{equation*}
and
\begin{equation*}
\prod_{i,j=1}^r(q^{1+m_j}z_i/z_j;q)_{-m_i}^{-1}=
\prod_{i,j=1}^r\frac{(qz_i/z_j;q)_{m_j}}{(qz_i/z_j;q)_{m_j-m_i}},
\end{equation*}
together with the $n\mapsto r$, $x_i\mapsto z_i$, and $y_i\mapsto-m_i$,
$i=1,\dots,r$, case of \cite[Lemma~3.12]{milne}, specifically
\begin{multline}
\prod_{i,j=1}^r(qz_i/z_j;q)_{m_j-m_i}=
(-1)^{(r-1)|{\mathbf m}|}
q^{-\binom{|\mathbf m|+1}2+r\sum_{i=1} ^r\binom{m_i+1}2}\\\times
\prod_{i=1}^rz_i^{|{\mathbf m}|-rm_i}
\prod_{1\le i<j\le r}\left(\frac {z_iq^{-m_i}-z_jq^{-m_j}}
{z_i-z_j}\right),
\end{multline}
and further the $a\mapsto aq^{-|{\mathbf m}|}$, $c\mapsto cq^{-|{\mathbf m}|}$,
$e_i\mapsto e_iq^{-m_i}$, $g\mapsto gq^{-|{\mathbf m}|}$,
and $z_i\mapsto z_iq^{-m_i}$, $i=1,\dots,r$,
case of the multidimensional summation formula in \eqref{p88gl},
we simplify the expression in \eqref{longlm} to
\begin{multline*}
(-1)^{|{\mathbf m}|}q^{-\binom{|\mathbf m|}2}
\left(\frac{fg}a\right)^{|{\mathbf m}|}
\prod_{i=1}^r\left(\frac{1-az_iq^{-m_i-|{\mathbf m}|}}{1-az_i}\right)\\\times
\prod_{i,j=1}^r\frac{(qz_i/z_j;q)_{m_j}}{(z_iq/e_iz_j;q)_{m_j}}
\prod_{i=1}^r\frac{(az_iq^{-m_i};q)_{-|{\mathbf m}|}\,(fz_i;q)_{-m_i}}
{(az_iq/e_i;q)_{-|{\mathbf m}|}\,(az_iq/g;q)_{-m_i}}\\\times
\frac{(g;q)_{-|{\mathbf m}|}}{(aq/f;q)_{-|{\mathbf m}|}}\,
\frac{(1-cq^{-|{\mathbf m}|-1})}{(1-c/q)}
\prod_{i=1}^r\frac{(1-az_iq^{1-m_i}/c)}{(1-az_iq/c)}\\\times
\left(1-\frac{(1-q^{|{\mathbf m}|}/g)(1-cf/aq)}{(1-c/gq)(1-Ef/a)}
\prod_{i=1}^r\frac{(1-ce_i/az_iq)}{(1-cq^{m_i}/az_iq)}\right)
\frac{(1-c/gq)}{(1-cq^{-|{\mathbf m}|-1})}\\\times
\frac{(a/Ef,aq/fg;q)_{\infty}}
{(aq^{|{\mathbf m}|}/Efg,aq^{1-|{\mathbf m}|}/f;q)_{\infty}}
\prod_{i=1}^r\frac{(az_iq^{1-m_i-|{\mathbf m}|},az_iq/e_ig;q)_{\infty}}
{(az_iq^{1-|{\mathbf m}|}/e_i,az_iq^{1-m_i}/g;q)_{\infty}}.
\end{multline*}
Now, this can easily be further transformed into
\begin{multline*}
\left(1-\frac{(1-q^{|{\mathbf m}|}/g)(1-cf/aq)}{(1-c/gq)(1-Ef/a)}
\prod_{i=1}^r\frac{(1-ce_i/az_iq)}{(1-cq^{m_i}/az_iq)}\right)
\prod_{i=1}^r\frac{(1-cq^{m_i}/az_iq)}{(1-c/az_iq)}\\\times
\frac{(1-c/gq)}{(1-c/q)}\frac{(a/Ef,aq/fg,q^{1+|{\mathbf m}|}/g;q)_{\infty}}
{(aq^{|{\mathbf m}|}/Efg,aq/f,q/g;q)_{\infty}}
\prod_{i,j=1}^r\frac{(qz_i/z_j,z_iq^{1+m_j}/e_iz_j;q)_{\infty}}
{(q^{1+m_j}z_i/z_j,z_iq/e_iz_j;q)_{\infty}}\\\times
\prod_{i=1}^r\frac{(az_iq,q/az_i,az_iq/e_ig,q^{1+m_i}/fz_i;q)_{\infty}}
{(az_iq/e_i,az_iq/g,q^{1+m_i}/az_i,q/fz_i;q)_{\infty}},
\end{multline*}
which is exactly the $b_i=aq^{-m_i}$, $i=1,\dots,r$, case of the right side
of \eqref{m88gl}.
\end{proof}

\end{document}